\documentclass{amsart}
\usepackage{amscd,amssymb,latexsym}

\theoremstyle{plain}
\newtheorem{thm}[subsection]{Theorem}
\newtheorem{lem}[subsection]{Lemma}
\newtheorem{prop}[subsection]{Proposition}
\newtheorem{cor}[subsection]{Corollary}

\theoremstyle{definition}
\newtheorem{rem}[subsection]{Remark}
\newtheorem{definition}[subsection]{Definition}
\newtheorem{example}[subsection]{Example}
\newtheorem{para}[section]{}

\newenvironment{alphenum}%
{%
 \begin{enumerate}
}%
{%
 \end{enumerate}%
}

\newenvironment{Alphenum}%
{%
 \begin{enumerate}
}%
{%
 \end{enumerate}%
}

\newenvironment{Alphenumprime}%
{%
 \begin{enumerate}
}%
{%
 \end{enumerate}%
}

\newenvironment{romenum}%
{%
 \begin{enumerate}
}%
{%
 \end{enumerate}%
}

\newenvironment{Romenum}%
{%
 \begin{enumerate}
}%
{%
 \end{enumerate}%
}

\newcommand{\abs}[1]{\left|#1\right|}
\newcommand{\tsum}{\textstyle\sum\limits}

\newcommand{\A}{\mathcal{A}}
\newcommand{\B}{\mathcal{B}}

\newcommand{\LL}{\mathcal{L}}
\newcommand{\KK}{\mathcal{K}}
\newcommand{\MM}{\mathcal{M}}
\newcommand{\RR}{\mathcal{R}}
\newcommand{\PP}{\mathcal{P}}

\newcommand{\C}{\mathbb{C}}
\newcommand{\N}{\mathbb{N}}
\newcommand{\F}{\mathbb{F}}
\newcommand{\Z}{\mathbb{Z}}
\newcommand{\K}{\mathbb{K}}
\newcommand{\R}{\mathbb{R}}

\renewcommand{\P}{\mathbb{P}}
\renewcommand{\SS}{\mathbb{S}}
\newcommand{\D}{\Delta}
\newcommand{\e}{\epsilon}
\newcommand{\g}{{\gamma}}

\newcommand{\s}{\sigma}

\newcommand{\bo}{\mathbf{0}}

\renewcommand{\a}{{\alpha }}
\renewcommand{\b}{{\beta }}
\renewcommand{\d}{\partial}
\renewcommand{\l}{\lambda}

\renewcommand{\L}{{\bigwedge}}

\DeclareMathOperator{\rank}{rank}
\DeclareMathOperator{\coker}{coker}
\DeclareMathOperator{\id}{id}

\DeclareMathOperator{\Hom}{Hom}

\DeclareMathOperator{\lk}{lk}

\DeclareMathOperator{\Tors}{Tors}
\DeclareMathOperator{\sgn}{sgn}

\DeclareMathOperator{\codim}{codim}

\DeclareMathOperator{\Mat}{Mat}

\begin{document}

\title[Cohomology rings and nilpotent quotients of arrangements]%
{Cohomology rings and nilpotent quotients of real and complex arrangements}

\author{Daniel~Matei}
\address{Department of Mathematics,
Northeastern University,
Boston, MA 02115}
\email{dmatei@lynx.neu.edu}

\author{Alexander~I.~Suciu}
\address{Department of Mathematics,
Northeastern University,
Boston, MA 02115}
\email{alexsuciu@neu.edu}
\urladdr{http://www.math.neu.edu/\~{}suciu}



\thanks{To appear in {\em Singularities and Arrangements, Sapporo--Tokyo 1998}, 
Advanced Studies in Pure Mathematics, Kinokuniya, Tokyo.  Available at 
{\ttfamily http://xxx.lanl.gov/abs/math.GT/9812087.}}

\dedicatory{Dedicated to Peter Orlik on his $60^{\text th}$ birthday}

\begin{abstract}
For an arrangement with complement $X$ and fundamental group $G$,
we relate the truncated cohomology ring, $H^{\le 2}(X)$, to the
second nilpotent quotient, $G/G_3$.  We define invariants of $G/G_3$
by counting normal subgroups of a fixed prime index $p$, according
to their abelianization.   We show how to compute this distribution
from the resonance varieties of the Orlik-Solomon algebra mod~$p$.
As an application, we establish the cohomology classification of
$2$-arrangements of $n\le 6$ planes in $\R^4$.
\end{abstract}

\maketitle

\section*{Introduction}
\label{sec:intro}

\para
Two basic homotopy-type invariants of a path-connected space $X$
are: the cohomology ring, $H^*(X)$, and the fundamental group, $G=\pi_1(X)$.
Given $X$ and $X'$, one would like to know:
\begin{Romenum}
\item \label{q1}
Is there an isomorphism $H^{\le{q}}(X)\cong H^{\le{q}}(X')$ between
the cohomology rings, up to degree $q$?
\item \label{q2}
Is there an isomorphism $G/G_{q+1}\cong G'/G'_{q+1}$ between the
$q^{\text{th}}$ nilpotent quotients?
\end{Romenum}

We single out a class of spaces---including complements of complex
hyperplane arrangements, complements of `rigid' links, and complements
of arrangements of transverse planes in $\R^4$---for which the above
questions are amenable to a detailed study, capable of yielding
classification results. The invariants that we use have a dual nature,
being able to capture both the  ring-theoretic properties of the cohomology
of $X$, and the group-theoretic properties of the nilpotent quotients of $G$.
Our main result is an explicit correspondence between two sets of
invariants---one determined by the vanishing cup products in $H^{\le 2}(X)$,
the other by the finite-index subgroups of $G/G_3$.

\para
For $q=1$, questions (\ref{q1}) and (\ref{q2}) are equivalent,
provided $H_1$ is torsion free.  Indeed, $H^1(X)=G/G_2$ under that assumption.
For $q=2$, the two questions are still equivalent, under some additional
conditions:  If $H_2$ is also torsion-free, and the cup-product map
$\mu:H^1\wedge H^1\to H^2$ is surjective, then:
\[
H^{\le{2}}(X)\cong H^{\le{2}}(X') \:\text{ if and only if }\:
G/G_{3} \cong G'/G'_{3}.
\]

Section~\ref{sec:cohonilp} is devoted to a proof of this fact.
A key ingredient is the vanishing of the Hurewicz map $\pi_2(X)\to H_2(X)$,
which permits us to identify $H^{\le 2}(X)$ with $H^{\le 2}(G)$.  The
other ingredient is the interpretation of the
$k$-invariant of the extension $0\to G_2/G_3\to G/G_3 \to G/G_2\to 0$,
in terms of the cup-products of $G$.

In Section~\ref{sec:pres}, we use commutator calculus
to describe the nilpotent quotients of $G$.
We restrict our attention to spaces $X$, for which $G=\pi_1(X)$
has a finite presentation $G=\F/R$, with $R\subset [\F,\F]$.
The cup products in $H^{\le{2}}(G)$ can then be computed
from the second order Fox derivatives of the relators.

\para
The invariants of the cohomology ring that we use are
the resonance varieties, first introduced
by Falk~\cite{Fa} in the context of complex hyperplane arrangements.
The $d^{\text{th}}$ resonance variety of $X$, with coefficients
in a field $\K$, is the set $\RR_d(X,\K)$ of cohomology classes
$\l \in H^1(X,\K)$ for which there is a subspace
$W \subset H^1(X,\K)$, of dimension $d+1$, such that $\mu(\l \wedge W)= 0$.

In Section~\ref{sec:resvar}, we prove that $\RR_d(X,\K)$ equals $\RR_d(G,\K)$,
the resonance variety of the Eilenberg-MacLane space $K(G,1)$.  Moreover,
we exploit the Fox calculus interpretation of cup products to show that
the varieties $\RR_d(G,\K)$ are the  determinantal varieties of
the `linearized' Alexander matrix of $G$.

\para
A well-known invariant of a group $G$ is the number of normal
subgroups of fixed prime index.  For a
commutator-relators group, that number depends only
on $n=\rank G/G_2$, and the prime $p$.  In order to
get a finer invariant, we consider the distribution of
index~$p$ subgroups, according to their abelianization.
The $\nu$-invariants of the nilpotent quotients $G/G_{q+1}$ are
defined as follows:
\[
\nu_{p,d}(G/G_{q+1}) = \#
\left\{ K\lhd G/G_{q+1}  \left|\:
\begin{array}
[c]{l}%
[G/G_{q+1}:K]=p \ \text{ and }\\
\dim_{\,\Z_p} (\Tors H_1(K))\otimes \Z_p = d
\end{array}
\right\}\right. .
\]

In Section~\ref{sec:pindex}, we show how to compute the
$\nu$-invariants of $G/G_{3}$ from the stratification of $\P(\Z_p^n)$
by the projectivized $\Z_p$-resonance varieties of $X$:
\[
\nu_{p,d}(G/G_3)=\#(\PP_d(X,\Z_p)\setminus\PP_{d+1}(X,\Z_p)).
\]
This formula makes the computation of the $\nu$-invariants practical.
It also makes clear that the mod $p$ resonance varieties of $X$,
which are defined solely in terms of $H^{\le 2}(X)$, do capture significant
group-theoretic information about $G/G_{3}$.

\para
In the case where $X$ is the complement of a complex hyperplane arrangement,
the varieties $\RR_d(X,\C)$ have been extensively studied by Falk, Yuzvinsky,
Libgober, Cohen, and Suciu \cite{Fa, Yuz, Li2, LY, CScv}.
The top variety, $\RR_1(X,\C)$, is a complete invariant of the
cohomology ring $H^{\le 2}(X)$.  Moreover, $\RR_1(X,\C)$ is a
union of linear subspaces intersecting only at the origin,
and $\RR_d(X,\C)$ is the union of those subspaces of dimension
at least $d+1$.

In Section~\ref{sec:cxa}, we use these results to derive
a simple consequence.  Since $\RR_d(X,\C)$ has integral equations,
we may consider its reduction mod $p$.  If that variety coincides
with $\RR_d(X,\Z_p)$, we have:
\[
\nu_{p,d-1}(G/G_3)=\frac{p^{d}-1}{p-1} m_{d},
\]
where $m_d$ is the number of components of $\RR_1(X,\C)$ of dimension $d$.

In general, though, this formula fails, due to a different resonance
at `exceptional' primes.  For such primes $p$, the variety
$\RR_d(X,\Z_p)$ is not necessarily the union of the components of
$\RR_1(X,\Z_p)$ of dimension $\ge d+1$.  Furthermore,
$\RR_1(X,\C)$ mod~$p$ and $\RR_1(X,\Z_p)$ may differ in
the number of non-local components, as well as in the
dimensions of those components. Most strikingly, $\RR_1(X,\Z_p)$
may have non-local components, even though $\RR_1(X,\C)$ mod~$p$ has none.

\para
Much of the original motivation for this paper came from an effort
to understand Ziegler's pair of arrangements of $4$ transverse planes
in $\R^4$. Those arrangements  have isomorphic lattices, but their
complements have non-isomorphic cohomology rings, see \cite{Zi}.
In an earlier work \cite{MS}, we investigated the homotopy types
of complements of $2$-arrangements, obtaining a complete classification
for $n \le 6$ planes.  This left open the problem of the classifying
cohomology rings for $n>4$.

In Section~\ref{sec:real}, we start a study of the
varieties $\RR_d(X,\C)$, where $X$ is the complement
of a $2$-arrangement. The resonance varieties of real arrangements
exhibit a much richer geometry than those of complex arrangements.
Most strikingly, $\RR_1(X,\C)$ may not be a union of
linear subspaces, and it may not determine $H^*(X)$.

Using the $\nu$-invariants of $G/G_3$, we establish the cohomology
classification of complements of $2$-arrangements of $n\le 6$ planes
in $\R^4$.  With one exception, this classification coincides with
the homotopy-type classification from \cite{MS}.  The exception
is Mazurovski\u{\i}'s pair \cite{Mz}. The two complements, $X$ and $X'$,
have isomorphic cohomology rings, and thus $G/G_3\cong G'/G'_3$.
On the other hand, $\nu_{3,2}(G/G_4)=162$ and $\nu_{3,2}(G'/G'_4)=172$.

As this example shows, the $\nu$-invariants of the third nilpotent
quotient {\em cannot} be computed from the resonance varieties of
the cohomology ring.  To arrive at a  more conceptual understanding of
these invariants, one needs to look beyond cup-products, and on
to the Massey products.  This will be pursued elsewhere.

\section*{Acknowledgment}
\label{sec:ack}
We wish to thank Sergey Yuzvinsky for valuable discussions
regarding the material in Section~\ref{sec:cxa}, and the 
referee for carefully reading the manuscript. 
The computations for this work were done with the help of the
packages {\it GAP~3.4.4} (\cite{Sch95}), {\it Macaulay~2} (\cite{GS}),
and {\it Mathematica~3.0}.

\setcounter{section}{0}
\section{Cohomology ring and second nilpotent quotient}
\label{sec:cohonilp}

In this section, we introduce a class of spaces that abstract
the cohomological essence of hyperplane arrangements.   We then
relate the cohomology ring of such a space $X$ to the
second nilpotent quotient of the fundamental group of $X$.

\subsection{Cohomology ring} All the spaces considered in this paper
have the homotopy type of a connected CW-complex with finite $2$-skeleton.
Let $X$ be such a space. Consider the following  conditions on the
cohomology ring of $X$:

\begin{Alphenum}
\item  \label{cond1}
The homology groups $H_*(X)$ are free abelian.
\item  \label{cond2}
The cup-product map $\mu_{X}:\L^* H^1(X)\to H^*(X)$
is surjective.
\end{Alphenum}
If conditions~(\ref{cond1}) and (\ref{cond2})
only hold for $1\le *\le n$, we will refer to them as
($\textup{\ref{cond1}}_n$) and ($\textup{\ref{cond2}}_n$).

\begin{example}
\label{exm:OSalgebra}
The basic example we have in mind is that of the complement,
$X=\C^{\ell}\setminus \bigcup_{H\in \A} H$, of a central hyperplane
arrangement $\A$ in $\C^{\ell}$.  As shown by Brieskorn~\cite{Bri}
(solving a conjecture of Arnol'd), such a space $X$ satisfies
conditions~(\ref{cond1}) and (\ref{cond2}).  Moreover, as shown by Orlik and
Solomon~\cite{OS}, the {\em intersection lattice} of the arrangement,
$L(\A)=\left\{\bigcap_{H\in \B}\, H\left| \B\subseteq\A\right\}\right.$,
determines the cohomology ring of $X$, as follows:
\[
H^{*}(X) = \left.\sideset{}{^*}\bigwedge \Z^n \right\slash \Big(
\partial e_{\B}\:\big|\,\codim \bigcap_{H\in \B} H < \abs{\B}\Big).
\]
Here $\bigwedge^{*}\Z^n$ is the exterior algebra on generators
$e_1,\dots ,e_n$ dual to the meridians of the hyperplanes; and, if
$\B=\{H_{i_1},\dots,H_{i_r}\}$ is a sub-arrangement, then
$e_{\B}=e_{i_1}\cdots e_{i_r}$, and
$\partial e_{\B}=\sum_{q}(-1)^{q}e_{i_1}\cdots\widehat{e_{i_q}}\cdots
e_{i_r}$.
See \cite{OT} for a thorough treatment of hyperplane arrangements.
\end{example}

Let $X$ be a space satisfying conditions ($\textup{\ref{cond1}}_n$) and
($\textup{\ref{cond2}}_n$). The first condition and the Universal Coefficient
Theorem (see~\cite{Bre}, Theorem~7.1, p.~281) imply that $H_*(X)=H^*(X)$, 
for $*\le n$. Write $H=H_1(X)=H^1(X)$. Denote by $I^*$ the kernel 
of the cup-product map. Condition~($\textup{\ref{cond2}}_n$)
can be restated as saying that the following sequence is exact:
\begin{equation}
\label{eq:cohopres}
0\to I^* \xrightarrow{\iota} \L^* H \xrightarrow{\mu_{X}} H^*(X) \to 0,
\quad \text{ for } *\le n.
\end{equation}
By condition~($\textup{\ref{cond1}}_n$), this is in fact a split-exact
sequence.

\subsection{Hurewicz homomorphism}
\label{subsec:hurewicz}
The following lemma was proved by Randell \cite{Ra} in the case where
$X$ is the complement of a complex hyperplane arrangement.

\begin{lem}
\label{lem:randell}
If $X$ satisfies conditions \textup{(\ref{cond1}$_n$)} and
\textup{(\ref{cond2}$_n$)}, then the Hurewicz homomorphism,
$h:\pi_{i}(X)\to H_{i}(X)$, is the zero map, for $2\le i\le n$.
\end{lem}

\begin{proof}  The proof is exactly as in \cite{Ra}.  Let
$p:\widetilde{X}\to X$ be the universal covering map. Recall
that $p_*:\pi_i(\widetilde{X})\to \pi_{i}(X)$ is an isomorphism, for $i\ge 2$.
By naturality of the Hurewicz map, universal coefficients,
and condition~(\ref{cond1}$_n$), it is enough to show that
$p^*:H^{i}(X)\to H^{i}(\widetilde{X})$ is the zero map.
This follows from $H^{1}(\widetilde{X})=0$, condition~(\ref{cond2}$_n$),
and the naturality of cup products:
$p^*\circ\mu_{X}=\mu_{\widetilde{X}}\circ \wedge^{i} p^*$.
\end{proof}

\subsection{Group cohomology}
\label{subsec:gpcoho}
Let $G$ be a group.  The (co)homology groups of $G$ are
by definition those of the corresponding Eilenberg-MacLane space:
$H_*(G)=H_*(K(G,1))$ and $H^*(G)=H^*(K(G,1))$.
Consider the following homological conditions on $G$:

\begin{Alphenumprime}
\item \label{condG1}
The homology groups $H_1(G)$ and $H_2(G)$ are finitely generated free abelian.

\item \label{condG2}
The cup-product map $\mu_{G}:H^1(G)\wedge H^1(G)\to H^2(G)$
is surjective.
\end{Alphenumprime}

\begin{prop}
\label{prop:xandg}
Let $X$ be a space satisfying conditions \textup{(\ref{cond1}$_2$)} and
\textup{(\ref{cond2}$_2$)}, and let $G=\pi_1(X)$ be its fundamental group.
Then the following hold:
\begin{alphenum}
\item \label{kga}
$H_1(G)\cong H_1(X)$ and $H_2(G)\cong H_2(X)$.
\item  \label{kgb}
The rings $H^{\le 2}(G)$ and $H^{\le 2}(X)$ are isomorphic.
\end{alphenum}
Therefore, $G$ satisfies conditions~\textup{(\ref{condG1})}
and~\textup{(\ref{condG2})}.
\end{prop}
\begin{proof}
Recall $X$ has the homotopy type of a connected $CW$-complex
$Y$ with finite $2$-skeleton. A $K(G,1)$ space may be obtained
from $Y$ by attaching suitable cells of dimension $\ge 3$.
The resulting map, $j:X\to K(G,1)$, induces an isomorphism
$H_1(X)\cong H_1(G)$.  From the Hopf exact sequence
$\pi_2(X)\xrightarrow{h} H_2(X)\to H_2(G)\to 0$ and
Lemma~\ref{lem:randell}, we get $H_2(X)\cong H_2(G)$.
This finishes the proof of (\ref{kga}).

By universal coefficients, the map $j^*: H^{i}(G)\to H^{i}(X)$
is a group isomorphism, for $i\le 2$.  By naturality of cup products,
we have $j^*\mu_{G}(a\wedge b)=\mu_{X}(j^*a\wedge j^*b)$.
This proves (\ref{kgb}).
\end{proof}

\begin{rem}
\label{rem:holonomy}
The above conditions on $X$ also appear in \cite{BP, Ry2}.
The surjectivity of $\mu:H^1(X)\wedge H^1(X)\to H^2(X)$
is stated there dually, as the injectivity of the holonomy map,
$\mu^{\top}:H_2(X)\to \L^2 H_1(X)$.
\end{rem}

\subsection{Nilpotent quotients}
\label{subsec:nilp}
Let $G$ be a finitely generated group.
The lower central series of $G$ is defined inductively by
$G_1=G$, $G_{q+1}=[G,G_q]$,
where $[G,G_q]$ denotes the subgroup of $G$ generated by the
commutators $[x,y]=xyx^{-1}y^{-1}$ with $x\in G$ and $y\in G_q$.
The quotient $G_{q}/G_{q+1}$ is a finitely generated abelian group,
called the {\em $q^{\text th}$~lower central series quotient} of $G$.
The quotient $G/G_{q+1}$ is a nilpotent group, called the
{\em $q^{\text th}$~nilpotent quotient} of $G$. See~\cite{MKS} for details.

We will be mainly interested in the second nilpotent quotient,
$G/G_3$.  This group is a central extension of finitely generated abelian
groups,
\begin{equation}
\label{eq:centralext}
0\to G_2/G_3 \to G/G_3 \to G/G_2 \to 0.
\end{equation}
The extension is classified by the $k$-invariant,
$\bar{\chi} \in H^2(G/G_2;G_2/G_3)$.
The isomorphism type of $G/G_3$ is determined by
$G/G_2$, $G_2/G_3$, and $\bar\chi$, as follows.

Let $G$ and $G'$ be two groups.
Then $G/G_3 \cong G'/G'_3$ if and only if there
exist isomorphisms $\phi:G/G_2\to G'/G'_{2}$ and
$\psi:G_{2}/G_{3}\to G'_{2}/G'_{3}$ under which
the $k$-invariants correspond:
$\psi_*(\bar\chi)=\phi^*(\bar\chi')\in H^2(G/G_2;G'_2/G'_3)$.

Now suppose $H=G/G_2$ is torsion-free.  As is well-known,
$H_*(H)\cong \L^* H$.  The {\em classifying map} for
the extension~\eqref{eq:centralext},
\[
\chi:\L^{2}H \to G_2/G_3
\]
is the image of $\bar\chi$ under the epimorphism
$H^2(H;G_2/G_3)\to \Hom(\L^2 H,G_2/G_3)$ provided by
the Universal Coefficient Theorem (see \cite{Dw}).
It is given by $\chi(x\wedge y) = [x,y]$ (see~\cite{Br},
Exercise~8, p.~97).
The condition that the $k$-invariants of $G/G_3$ and $G'/G'_3$
correspond translates to $\psi\circ\chi=\chi'\circ\wedge^{2}\phi$.
We shall write this equivalence relation between classifying maps
as $\chi\thicksim \chi'$.

Suppose now that $G_2/G_3$ is also torsion-free. Then,
the universal coefficient map is an isomorphism, and so
$\bar{\chi}$ and $\chi$ determine each other. Thus, for
a group $G$ with $G/G_2$ and $G_2/G_3$ torsion-free,
the isomorphism type of $G/G_3$ is completely determined
by the equivalence class of the classifying map $\chi$.

\subsection{Cup product and commutators}
\label{subsec:cupandcomm}
The $5$-term exact sequence for the extension
$0\to G_2 \to G \xrightarrow{\a} G/G_2 \to 0$ yields:
\begin{equation}
\label{eq:eta}
H_2(G) \xrightarrow{\a_*} H_2(G/G_2)
\xrightarrow{\delta} G_2/G_3 \to 0.
\end{equation}
Under the identification $H_2(G/G_2)\cong \L^2 H$, the
boundary map $\delta$ corresponds to the classifying map
$\chi$ (see \cite{Br}, Exercise~6, p.~47).
The next lemma interprets the map $\a_*$
in terms of the ring structure of $H^*(G)$.

\begin{lem}
\label{lem:etamu}
The map $\alpha_*:H_2(G) \to \L^2 H$ is the dual of the
cup-product map $\mu_G: H^1(G)\wedge H^1(G)\to H^2(G)$.
\end{lem}

\begin{proof}
Follows from the commutativity of the diagram
\begin{equation*}
\begin{CD}
& H^1(H)\wedge  H^1(H)   @>\a^*\wedge \a^*>> H^1(G)\wedge  H^1(G)\\
&  @VV{\mu_H}V               @VV{\mu_G}V \\
& H^2(H) @>\a^*>>    H^2(G)\\
\end{CD}
\end{equation*}
and the fact that the top and left arrows are isomorphisms.
\end{proof}

The following proposition generalizes a result proved by
Massey and Traldi~\cite{MT} in the case where $G$ is a link group.

\begin{prop}
\label{prop:masseytraldi}
Let $G$ be a group satisfying conditions \textup{(\ref{condG1})}
and \textup{(\ref{condG2})}.  Then $G_2/G_3$ is torsion-free,
and the following is a split exact sequence:
\begin{equation}
\label{eq:muchi}
0\to H_2(G) \xrightarrow{\mu^{\top}}  \L^2 H
\xrightarrow{\chi} G_2/G_3 \to 0.
\end{equation}
\end{prop}

\begin{proof}
The proof follows closely that in~\cite{MT}.
By Lemma~\ref{lem:etamu}, sequence~\eqref{eq:eta}
can be written as
$
H_2(G) \xrightarrow{\mu^{\top}}  \L^2 H
\xrightarrow{\chi} G_2/G_3 \to 0.
$
By condition~(\ref{condG2}), the map $\mu^{\top}$ is a
monomorphism, whence the exactness of \eqref{eq:muchi}.

Since $\mu: \L^2 H \to H^2(G)$ is an epimorphism between
finitely generated free abelian groups, it admits a splitting.
Hence $\mu^{\top}$ is a split injection, and so $\chi^{\top}$
is a split surjection.  Since $\L^2 H$ is torsion-free,
$G_2/G_3$ is also torsion-free.
\end{proof}

\begin{rem}
\label{rem:dwyer}
The injectivity of  $\mu^{\top}:H_2(G) \to  \L^2 H$ is
equivalent to the vanishing of $\Phi_3(G)$, where
$H_2(G)=\Phi_2(G)\supset \Phi_3(G) \supset\cdots$ is the Dwyer 
filtration, $\Phi_k(G) = \ker(H_2(G)\to H_2(G/G_{k-1}))$, 
see \cite{Dw}.  
\end{rem}

\subsection{Isomorphisms}
\label{subsec:isos}
The next result is an immediate
consequence of Proposition~\ref{prop:masseytraldi}:

\begin{prop}
\label{prop:iandg3}
Let $X$ be a space satisfying conditions \textup{(\ref{cond1}$_2$)} and
\textup{(\ref{cond2}$_2$)}, and let $G=\pi_1(X)$.
Then $I^2=G_2/G_3$, and the exact sequence
\begin{equation}
\label{eq:chimu}
0\to I^2 \xrightarrow{\iota} \L^2 H^1(X)
\xrightarrow{\mu} H^2(X) \to 0
\end{equation}
is the dual of sequence \eqref{eq:muchi}.
\end{prop}

We are now ready to establish the correspondence between
the truncated cohomology ring of $X$ and the second nilpotent
quotient of $G=\pi_1(X)$.  A version of the
equivalence $(\ref{eqb}) \Leftrightarrow (\ref{eqc})$ below,
with $\chi$ replaced by $\mu^{\top}$, was first established
by Traldi and Sakuma~\cite{TS}, in the case where
$X$ is a link complement.

\begin{thm}
\label{thm:traldisakuma}
Let $X$ and $X'$ be two spaces satisfying conditions
\textup{(\ref{cond1}$_2$)} and \textup{(\ref{cond2}$_2$)},
and let $G$ and $G'$ be their fundamental groups.
The following are equivalent:
\begin{alphenum}
\item	\label{eqa} $H^*(X) \cong H^*(X')$ for $*\le 2$;
\item \label{eqb} $G/G_3 \cong G'/G'_3$;
\item \label{eqc} $\chi \thicksim \chi'$.
\end{alphenum}
\end{thm}

\begin{proof}
$(\ref{eqa}) \Leftrightarrow (\ref{eqc})$.
By Proposition~\ref{prop:iandg3}, sequence ~(\ref{eq:chimu}) is exact, and
$\chi=\iota^{\top}$. The equivalence follows from the definitions.

$(\ref{eqb}) \Leftrightarrow (\ref{eqc})$.
By Propositions~\ref{prop:xandg} and \ref{prop:iandg3},
the first two lower central series quotients of $G$ and $G'$
are torsion-free. The equivalence follows from the discussion
in~\ref{subsec:nilp}.
\end{proof}

\subsection{Invariants of $H^{\le 2}(X)$ and $G/G_3$}
\label{subsec:handginv}
In view of Theorem~\ref{thm:traldisakuma}, an invariant of either
the truncated cohomology ring $H^{\le 2}(X)$, or the second nilpotent
quotient $G/G_3$, or the classifying map $\chi$, is an invariant of
the other two.  We will define in the latter sections a series of
invariants of both $H^{\le 2}(X)$  and $G/G_3$, and relate them one
to another.  For now, let us define invariants of $\chi$,
following an idea of Ziegler~\cite{Zi}, that originated
from Falk's work on minimal models of arrangements~\cite{Fa1}.

Let $\mu_{H}:\L^{i} H\otimes \L^{j} H\to \L^{i+j} H$ be the multiplication
in the exterior algebra $\L^* H$.
Consider the following finitely generated abelian group:
\[
Z_{i,j}(\chi)=\coker\bigg(\L^{i}H \otimes\L^{j}G_2/G_3
\xrightarrow{\id\otimes\L^{j}\chi^{\top}} \L^{i}H \otimes \L^{2j}H
\xrightarrow{\mu_{H}}\L^{i+2j}H\bigg).
\]
Clearly, if $\chi \thicksim \chi'$ then $Z_{i,j}(\chi) \cong Z_{i,j}(\chi')$.
Thus, the rank and elementary divisors of $Z_{i,j}(\chi)$ provide
invariants of both
$H^{\le 2}(X)$ and $G/G_3$.

\section{Generators and Relators}
\label{sec:pres}

In this section, we write down explicitly some of the maps introduced
in the previous section.  We start with a review of some basic facts
about Hall commutators and the Fox calculus.

\subsection{Basic commutators}
\label{subsec:hall}
Let $\F(n)$ be the free group on generators $x_1, \dots, x_n$.
A {\em basic commutator} in $\F=\F(n)$ is defined inductively as follows
(see \cite{Fe, MKS}):
\begin{alphenum}
\item \label{halla}
Each basic commutator $c$ has length $\ell(c)\in\N$.
\item \label{hallb}
The basic commutators of length $1$ are the generators $x_1,\dots,x_n$;
those of length $>1$ are of the form $c=[c_1,c_2]$, where $c_1$, $c_2$
are previously defined commutators and $\ell(c)=\ell(c_1)+\ell(c_2)$.
\item  \label{hallc}
Basic commutators of the same length are ordered arbitrarily;
if $\ell(c)>\ell(c')$, then $c>c'$.
\item  \label{halld}
If $l(c)>1$ and $c=[c_1,c_2]$, then $c_1<c_2$;
if $l(c)>2$ and $c=[c_1,[c_2,c_3]]$, then $c_1\ge c_2$.
\end{alphenum}

The basic commutators of the form
$c=[x_{i_1},[x_{i_2},[\dots[x_{i_{q-1}},x_{i_q}]\dots]]]$
are called {\em simple}.  We shall write them as
$c=[x_{i_1}, x_{i_2},\dots, x_{i_q}]$.
For $q\le 3$, all basic commutators are simple.

The following theorem of Hall is well-known (see {\it loc. cit.}):

\begin{thm}
The group $\F_q/\F_{q+1}$ is free abelian, and has a basis
consisting of the basic commutators of length $q$.
\end{thm}

In particular, if $w\in \F$ and $c_1, \dots, c_r$ are the basic
commutators of length $< q$, then $w^{(q)}:=w \bmod \F_q$ may be written
uniquely as $w^{(q)}=c_1^{e_1}c_2^{e_2} \cdots c_r^{e_r}$, for some
integers $e_1, \dots , e_r$.

The Hall commutators may be used to write down presentations
for the nilpotent quotients of a finitely presented group $G=\F/R$.
Indeed, if $G=\langle x_1,\dots,x_n\mid r_1, \dots, r_m\rangle$,
we have the following presentation for $G/G_q=\F/R\F_q$:
\begin{equation}
\label{eq:gmodgq}
G/G_q=\langle x_1, \dots,x_n\mid r_1^{(q)},\dots, r_m^{(q)},
c_1,\dots,c_l \rangle,
\end{equation}
where $r_k^{(q)}=r_k\mod\F_q$, and $\{c_h\}_{1\le h\le l}$ are the
basic commutators of length $q$.

\subsection{Fox calculus}
\label{subsec:fox}
Let $\Z\F$ be the group ring of $\F$, with augmentation map
$\e:\Z\F \to \Z$ given by $\e(x_i)=1$.
To each $x_i$ there corresponds a {\em Fox derivative}, $\partial_i:\Z\F
\to \Z\F$,
given by $\partial_i(1)=0$, $\partial_i(x_j)=\delta_{ij}$ and
$\partial_i(uv)=\partial_i(u)+\e(u)\partial_i(v)$.
The higher Fox derivatives, $\partial_{i_1,\dots ,i_k}$, are
defined inductively in the obvious manner.  The composition
of the augmentation map with the higher derivatives yields operators
$\e_{i_1,\dots ,i_k}:=\e\circ\partial_{i_1,\dots,i_k}:\Z\F\to\Z$.

Let $\a:\F(n)\to\Z^n$ be the abelianization map, given by $\a(x_i)=t_i$.
The following lemma is left as an exercise in the definitions.

\begin{lem} \label{lem:foxcalc}
We have:
\begin{alphenum}
\item \label{foxa}
$\partial_i [u,v]=(1-uvu^{-1})\partial_i u+(u-[u,v])\partial_i v$.
\item \label{foxc}
$\a(\partial_i[x_{i_1}, x_{i_2},\dots, x_{i_q}])=(t_{i_1}-1)\cdots
(t_{i_{q-2}}-1)
\big((t_{i_{q-1}}-1)
\delta_{i,i_q}- (t_{i_q}-1) \delta_{i,i_{q-1}}\big).$
\item \label{foxd}
$\e_I(w)=0$, if $w\in\F_q$ and $|I|<q$.
\item \label{foxe}
$\e_I(uv)=\e_I(u)+\e_I(v)$, if $u,v\in\F_q$ and $|I|=q$.
\end{alphenum}
\end{lem}

\subsection{Commutator relations}
\label{subsec:gppres}
We now make more explicit some of the constructions from
section~\ref{sec:cohonilp}, for the following class of groups.

\begin{definition}
\label{def:commrel}
A group $G$ is called a {\em commutator-relators group} if it
admits a presentation $G=\F(n)/R$, where $R$ is the normal
closure of a finite subset of $[\F,\F]$.
\end{definition}
In other words, $G$ has a finite presentation
$G=\langle x_1,\dots,x_n\mid r_1,\dots ,r_m\rangle$, and $G/G_2=\Z^n$.
Commutator-relators groups appear as fundamental groups of certain
spaces that we shall encounter later on.  The following proposition
gives sufficient conditions for this to happen.

\begin{prop}
\label{prop:twocw}
Let $X$ be a space that is homotopy equivalent to a finite
CW-complex $Y$, with $1$-skeleton $Y^{(1)}=\bigvee_{i=1}^{n} S^1_i$.
If $H_1(X)=\Z^n$, then $G=\pi_1(X)$ is a commutator-relators group.
\end{prop}

\begin{proof}
The $2$-skeleton $Y^{(2)}=\bigvee_{i=1}^{n} S^1_i \cup\bigcup_{k=1}^{m}
e^2_k$
determines a presentation $G=\langle x_1,\dots ,x_n\mid r_1,\dots
,r_m\rangle$.
A presentation matrix for the abelianization of $G$ is $E=(\e_i(r_k))$.
Since $H_1(X)=\Z^n$, we have $H_1(G)=\Z^n$.  Thus, $E$ is equivalent to
the zero matrix, and hence $E$ {\em is} the zero matrix. Thus,
all relators $r_k$ are commutators.
\end{proof}

Now let $\phi:\F\to G$ be the quotient map, and let $\a:G\to G/G_2$ be the
abelianization map. Set $t_i=\a(\phi(x_i))$.  Then $\{t_1,\dots,t_n\}$
form a basis for $H_1(G)$, and their Kronecker duals,
$\{e_1,\dots ,e_n\}$, form a basis for $H^1(G)$.

By the Hopf formula, we have $H_2(G)=R/[R,\F]$.
Assume that $H_2(G)$ is free abelian, and let $\theta_k=r_k \bmod [R,\F]$.
Then $\{\theta_1,\dots,\theta_m\}$  form a basis for $H_2(G)$, and their
duals, $\{\g_1,\dots ,\g_m\}$, form a basis for $H^2(G)$.

\begin{prop}
\label{prop:cupprod}
Let $G$ be a commutator-relators group, with $H_2(G)$ free abelian.
In the basis specified above, the cup-product map
$\mu:H^1(G)\wedge H^1(G)\to H^2(G)$ is given by
\[\mu(e_i \wedge e_j)=\tsum_{k=1}^{m}{\e_{i,j}(r_k)\g_k}.\]
\end{prop}
\begin{proof}  This follows immediately from~\cite{FS}, Theorem~2.3.
\end{proof}

\subsection{Links in $\SS^3$}
\label{subsec:links}
We conclude this section with a classical example.
Let $L=\{L_1,\dots , L_n\}$ be an oriented link in $\SS^3$.
Its complement, $X=\SS^3\setminus \bigcup_i L_i$, has the homotopy
type of a connected, $2$-dimensional finite CW-complex.  The homology
groups of $X$ are computed by Alexander duality:  $H_1(X)=\Z^n$,
$H_2(X)=\Z^{n-1}$. It follows that Condition~(\ref{cond1}) is always
satisfied for a link complement.  If $L=\hat\beta$ is the closure of a pure
braid
$\beta\in P_n$, then $X$ satisfies the assumption of
Proposition~\ref{prop:twocw},
and so $G=\pi_1(X)$ is a commutators-relators group, with presentation
$G=\langle x_1,\dots ,x_n \mid \beta(x_i)x_i^{-1}=1,\ 1\le i <n\rangle$.

For an arbitrary link $L$, let $\{e_1,\dots,e_n\}$ be the
basis for $H^1(X)$ dual to the meridians of $L$.
Choose arcs $c_{i,j}$ in $X$ connecting
$L_i$ to $L_{j}$, and let $\gamma_{i,j}\in H^2(X)$ be their duals.
Then $\{\gamma_{1,n},\dots ,\gamma_{n-1,n}\}$ forms a basis for
$H_2(X)$.  Let $l_{i,j}=\lk(L_i,L_j)$ be the linking
numbers of $L$. A presentation for the cohomology ring
of $X$ is given by:
\begin{equation}
\label{eq:coholink}
H^{*}(X) = \left( e_i, \gamma_{i,j} \: \bigg|\:
\begin{array}{l}%
e_i e_j = l_{i,j} \gamma_{i,j},\
\gamma_{i,j} + \gamma_{j,k} + \gamma_{k,i} = 0\\
\gamma_{i,i}=\gamma_{i,j}e_k=e_k\gamma_{i,j}=\gamma_{i,j}\gamma_{k,l}=0
\end{array}
\right).
\end{equation}

Let $\mathcal{G}$ be the ``linking graph" associated to $L$:
It is the complete graph on $n$ vertices, with edges labelled
by the linking numbers. If $\mathcal{G}$ possesses a spanning
tree $T$ with $n$ vertices, and all edges labelled $\pm 1$,
we say that $L$ is {\em (cohomologically) rigid}. The complement
of such a link satisfies condition~(\ref{cond2}), see~\cite{MT, La, BP}.
Moreover, $G_2/G_3$ is free abelian of rank $\binom{n-1}{2}$,
with basis $\{x_{ij} \mid ij \not\in T \text{ and } i<j\}$.  The classifying
map $\chi\colon\L^{2}H \to
G_2/G_3$ is given by
\[
\chi(e_i \wedge e_j)=
\begin{cases}
x_{ij} & \text{if }\, ij \not\in T, \\
\sum_{\{k\mid ik\not\in T\}}{l_{i,k}x_{ik}} & \text{if }\, ij \in T,
\end{cases}
\]
where $x_{ik}=-x_{ki}$, for $i>k$.

We will be mainly interested in those rigid links for which $l_{i,j}=\pm 1$.
Examples  include the Hopf links, and, more generally, the singularity links
of $2$-arrangements in $\R^4$ (see \ref{subsec:planes}).
For such links, the presentation~\eqref{eq:coholink} simplifies to:
\begin{equation}
\label{eq:linkcohoring}
H^{*}(X) = \left( e_i \mid
e_i^2=0,\ e_ie_j=-e_je_i,\ l_{i,j}e_ie_j + l_{j,k}e_je_k + l_{k,i}e_ke_i =
0\right).
\end{equation}
Moreover, the transpose of the classifying map,
$\chi^{\top}: G_2/G_3 \to \L^{2}H$, is given by the simple formula
\begin{equation}
\label{eq:chitop}
\chi^{\top}(x_{ij})=(e_{i}-l_{i,j}e_{n})\wedge (e_j-l_{i,j}e_{n}).
\end{equation}

\section{Resonance varieties}
\label{sec:resvar}

In this section, we define the `resonance' varieties of
the cohomology ring of a space $X$.  We then show that,
under certain conditions on $X$, these varieties are the
determinantal varieties of the linearized Alexander matrix
of the group $G=\pi_1(X)$.

\subsection{Filtration of first cohomology}
\label{subsec:resfil}
Let $X$ be a space that satisfies conditions
($\textup{\ref{cond1}}_2$) and ($\textup{\ref{cond2}}_2$)
of Section~\ref{sec:cohonilp}, and the hypothesis of
Proposition~\ref{prop:twocw}.  We thus have: $H^1(X)=\Z^n$,
$H^2(X)=\Z^m$, the cup-product map $\mu:H^1(X)\wedge H^1(X)\to H^2(X)$
is surjective, and $G=\pi_1(X)$ is a commutator-relators group.

\begin{lem}
\label{lem:uct}
Let $X$ be as above, and let $\K$ be a commutative field. 
\begin{alphenum}
\item \label{Kcup}
The $\K$-cup products may be computed from the integral ones:
$\mu_{\K}=\mu\otimes\id_{\K}$.
\item \label{ZtoK}
If $H^{\le 2}(X)\cong H^{\le 2}(X')$ then
$H^{\le 2}(X,\K)\cong H^{\le 2}(X',\K)$.
\end{alphenum}
\end{lem}

\begin{proof}
Let $\kappa: \Z\to \K$ be the homomorphism given by $\kappa (1)=1$. 
From the definitions, the coefficient map 
$\kappa_*: H^{*}(X, \Z)\to H^{*}(X, \K)$, and the map 
$\id\otimes\,\kappa: H^{*}(X)\otimes\Z\to H^{*}(X)\otimes \K$  
commute with cup products.  By the Universal Coefficient Theorem 
(see \cite{Bre}, Theorem~7.4, p.~282), the map 
$\upsilon: H^{*}(X)\otimes\K\to H^{*}(X,\K)$, 
$\upsilon([z]\otimes k)=[z\otimes k]$ is an isomorphism 
for $*\le 2$.  Since $\upsilon\circ (\id\otimes\,\kappa)=\kappa_*$, 
the map $\upsilon$ also commutes with cup products.  
The conclusions follow.
\end{proof}

\begin{definition}
\label{def:resvar}
Let $d$ be an integer, $0\le d \le n$.
The $d^{\text{th}}$~{\em resonance variety}~of $X$ (with coefficients
in $\K$) is the subvariety of $H^1(X,\K)=\K^n$, defined as follows:
\[
\RR_{d}(X,\K) =\left\{ \l \in H^1(X,\K)\: \left|\:
\begin{array}
[c]{l}%
\exists \text{ subspace $W \subset H^1(X,\K)$ such that}\\
 \dim W=d+1 \text{ and } \mu(\l \wedge W)= 0
\end{array}
\right\}\right. .
\]
\end{definition}

The resonance varieties form a descending filtration
$\K^n=\RR_0\supset \RR_1\supset\cdots\supset \RR_{n-1}\supset
\RR_n=\emptyset$.
The {\em ambient type} of the $\K$-resonance varieties depends only on the
truncated cohomology ring $H^{\le 2}(X,\K)$, and thus, 
by Lemma~\ref{lem:uct}~(\ref{ZtoK}), only on $H^{\le 2}(X)$. More precisely, if
$H^{\le 2}(X) \cong H^{\le 2}(X')$, there exists a linear automorphism
of $\K^n$ taking $\RR_{d}(X,\K)$ to $\RR_{d}(X',\K)$.

For a group $G$, define the resonance varieties to be those
of the corresponding Eilenberg-MacLane space:
$\RR_{d}(G,\K):=\RR_{d}(K(G,1),\K)$.

\begin{prop}
\label{prop:reskg1}
Let $X$ be a space satisfying conditions
\textup{(\ref{cond1}$_2$)} and \textup{(\ref{cond2}$_2$)}.
Let $G=\pi_1(X)$.  Then $\RR_{d}(X,\K)=\RR_{d}(G,\K)$.
\end{prop}

\begin{proof}
By Proposition~\ref{prop:xandg},
the inclusion $j:X\to K(G,1)$ induces an isomorphism
$j^{*}: H^{\le 2}(G)\to H^{\le 2}(X)$.  The conclusion
follows from remark (\ref{ZtoK}) above.
\end{proof}

\subsection{Alexander matrices}
\label{subsec:linalex}
Let $G=\langle x_1,\dots ,x_n\mid r_1,\dots ,r_m\rangle$ be a
commutator-relators group.  Recall the projection map
$\phi:\F(n)\to G$,  and the abelianization map, $\a:G\to \Z^n$,
given by $\a(x_i)=t_i$.

\begin{definition}
\label{def:alexmat}
The {\em Alexander matrix}	of $G$ is
the $m\times n$ matrix $A=(\a\phi\partial_i(r_k))$ with entries in the
Laurent polynomials ring $\Z[t_1^{\pm 1}, \dots,t_n^{\pm 1}]$.
\end{definition}

Now let $\psi:\Z[t_1^{\pm 1}, \dots, t_n^{\pm 1}] \to
\Z[[s_1,\dots, s_n]]$ be the ring homomorphism given by $\psi(t_i)=1+s_i$
and $\psi(t_i^{-1})=\sum_{q\ge 0}(-1)^{q} s_{i}^{q}$.  Also, let $\psi^{(q)}$
be the graded $q^{\text{th}}$ piece of $\psi$. Since all the relators of $G$
are commutators, the entries of $A$ are in the ideal $(t_1-1,\dots,
t_n-1)$, and so
$\psi^{(0)}A$ is the zero matrix.

\begin{definition}
\label{def:linalex}
The {\em linearized Alexander matrix} of $G$ is the
$m\times n$ matrix
\[
M=\psi^{(1)}A .
\]
\end{definition}

Note that the entries of $M$ are integral linear forms in $s_1,\dots, s_n$.
By Lemma~\ref{lem:foxcalc} (\ref{foxa}),~(\ref{foxc})  we have
$\psi^{(1)}\a\phi\partial_i(r_k)=\psi^{(1)}\a\phi\partial_i(r_k^{(3)})$.
Thus, $M$ depends only on the relators of $G$, modulo length~$3$ commutators.
By Lemma~\ref{lem:foxcalc}~(\ref{foxd}), (\ref{foxe}) those truncated relators
are given by $r_k^{(3)}=\prod_{i<j}[x_i,x_j]^{\e_{i,j}(r_k)}$.
Thus, the entries of $M$ are:
\begin{equation}
\label{eq:linalex}
M_{k,j}=\sum_{i=1}^{n}\e_{i,j}(r_k) s_i.
\end{equation}

The linearized Alexander matrix of a link was first considered
by Traldi~\cite{Tr}.  If the link $L$ has $n$ components,
then $M$ has size $n\times (n-1)$, and its entries are
$M_{k,j}=l_{k,j}s_k-\delta_{k,j}(\sum_{i}l_{k,i}s_i)$.

\subsection{Equations for resonance varieties}
\label{subsec:eqrv}
We now find explicit equations for the varieties $\RR_d(X,\K)$.
In view of Proposition~\ref{prop:reskg1}, that is the same
as finding equations for $\RR_d(G,\K)$, with $G=\pi_1(X)$.
Moreover, in view of Lemma~\ref{lem:uct}~(\ref{Kcup}), the formula
for the integral cup products from Proposition~\ref{prop:cupprod}
may be used to compute the $\K$-cup products. We will use the
notations of that proposition for the rest of this section.

Let $M$ be the linearized Alexander matrix of $G$.
Let $M_{\K}$ be the corresponding matrix of linear
forms over $\K$, and let $M(\l)$ be the matrix
$M_{\K}$ evaluated at $\l=(\l_1,\dots,\l_n)\in\K^n$.

\begin{thm}
\label{thm:reseqts}
For $G$ a commutator-relators group with $H_2(G)$ torsion free,
\[
\RR_{d}(G,\K)=\{\l\in \K^n \mid \rank_{\K} M(\l) < n-d\: \}.
\]
\end{thm}

\begin{proof}
Let $\l=\sum_{i=1}^n \l_ie_i\in H^1(G,\K)=\K^n$.
We are looking for $v=\sum_{i=1}^n v_ie_i$
such that $\mu(\l\wedge v)=0$ in $H^2(G,\K)=\K^m$.
Recall from Proposition~\ref{prop:cupprod} that
$\mu(e_i\wedge e_j)=\sum_{k=1}^{s}{\e_{i,j}(r_k)\g_k}$.
It follows that
\[
\mu(\l\wedge v)=\sum_{k=1}^{s} \bigg( \sum_{1\le i,j\le n}
\l_i v_j \e_{i,j}(r_k) \bigg) \g_k.
\]
We thus obtain a linear system of $s$ equations in $v_1,\dots, v_n$:
\begin{equation*}
\label{eq:Mmatrix}
\sum_{j=1}^{n}\bigg( \sum_{i=1}^{n}\l_i\e_{i,j}(r_k) \bigg) v_j=0,
\end{equation*}
with coefficients matrix $M(\l)$.

Now $\l$ belongs to $\RR_{d}(G,\K)$
if and only if the space $W$ of solutions of the linear system
$M(\l)\cdot v=0$ is at least $(d+1)$-dimensional. That translates
into the condition $\rank_{\K} M(\l) < n-d$ of the statement, and we are
done.
\end{proof}

We will be mainly interested in the coefficients fields
$\K=\C$ and $\K=\Z_p$, for some prime $p$.  By the above
theorem, the $\C$-resonance varieties  have {\em integral}
equations.  As we shall see in Section~\ref{sec:cxa}, although
$\mu_{\Z_p}:H^1(X,\Z_p)\wedge H^1(X,\Z_p)\to H^2(X,\Z_p)$ is
the reduction $\bmod\,p$ of $\mu:H^1(X)\wedge H^1(X)\to H^2(X)$,
the variety $\RR_{d}(X,\Z_p)$ is {\em not} necessarily the reduction
$\bmod\,p$ of $\RR_{d}(X,\C)$.

\begin{example}
\label{ex:rvlink}
Let $X$ be the complement of an $n$-component rigid link.
The matrix $M(\l)$ has entries
$M(\l)_{k,j}=l_{k,j}\l_k-\delta_{k,j}(\sum_{i}l_{k,i}\l_i)$.
The variety $\RR_1(X,\K)$ is the zero-locus of a degree~$n-2$
homogeneous polynomial obtained by taking the greatest common
divisor of the $(n-1)\times (n-1)$ minors of the matrix $M(\l)$.
At the other extreme, we have $\RR_{n-1}(X,\K)=\{\bo\}$.
Indeed, the off-diagonal entries of $M(\l)$ corresponding
to the edges of the maximal spanning tree generate the
maximal ideal $(\l_1,\dots ,\l_n)$ of $\K[\l_1,\dots,\l_n]$.
\end{example}

\subsection{Projectivized resonance varieties}
\label{subsec:projres}
The affine variety $\RR_{d}(X,\K)\subset \K^n$ is homogeneous,
and so defines a projective variety $\PP_{d}(X,\K)\subset \P(\K^n)$.
If $H^{\le 2}(X)$ is isomorphic to $H^{\le 2}(X')$, there is a
projective automorphism $\P(\K^n)\to\P(\K^n)$ taking $\PP_{d}(X,\K)$
to $\PP_{d}(X',\K)$.  The rest of the above discussion applies
to the projective resonance varieties in an obvious manner.
In particular, we have:

\begin{cor}
\label{cor:projres}
$
\PP_{d}(G,\K)=\{\l\in \P(\K^n) \mid \rank_{\K} M(\l) < n-d-1\}.
$
\end{cor}

\section{Prime index normal subgroups}
\label{sec:pindex}
In this section, we consider nilpotent quotients
of commutator-relators groups. We show how to count
the normal subgroups of prime index, according to
their abelianization.

\subsection{Counting subgroups}
\label{subsec:count}
Let $G$ be a group.  For a prime number $p$, let $\Sigma_p(G)$ be
the set of index~$p$ normal subgroups of $G$, and let
$N_p(G)=|\Sigma_p(G)|$ be its cardinality.

\begin{prop}
\label{prop:npfn}
For the free group $\F(n)$, the set
$\Sigma_p(\F(n))$ is in bijective correspondence with
the projective space $\P(\Z_p^n)$.
\end{prop}

\begin{proof}
Every index $p$ normal subgroup of $\F_n$ is the kernel of
an epimorphism $\lambda:\F(n)\to \Z_p$.  Such
homomorphisms are parametrized by $\Z_p^n \setminus \{\bo\}$.
Two epimorphisms $\lambda$ and $\lambda'$ have the
same kernel if and only if $\lambda=q\cdot \lambda'$, for
some $q\in \Z_p^*$.
\end{proof}

\begin{cor}
\label{cor:number}
Let $G=\F(n)/R$ be a commutator-relators group. For all primes $p$,
\[
N_p(G)=\frac{p^{n}-1}{p-1}.
\]
\end{cor}

\begin{proof}
Since $R$ consists of commutators, $\Hom(G,\Z_p)\cong\Hom(\F(n),\Z_p)$.
Thus, $\Sigma_p(G)$ is in one-to-one correspondence with
$\Sigma_p(\F(n))=\P(\Z_p^n)$.
\end{proof}

\subsection{Abelianizing normal subgroups}
\label{subsec:abalg}
Let $G=\langle x_1,\dots ,x_n\mid r_1,\dots ,r_m\rangle$ be
a commutator-relators group. Let $K\lhd G$ be a normal subgroup
of index~$p$, defined by a representation $\lambda:G\to \Z_p$,
$\l(x_i)=\l_i$.
Let $\bar\lambda:\Z G\to \Z\Z_p$ be the linear extension of $\l$
to group rings.  More precisely, we view here $\Z_p$ as a
multiplicative group, with generator $\zeta$.  Then
$\bar\l (x_i)=\zeta^{\l_i}$.
Finally, let $\b:\Z\Z_p \to \Mat(p,\Z)$ be the ring homomorphism
defined by the (left) regular representation of $\Z_p$.

\begin{definition}
\label{def:twistalex}
For a given representation $\l:G\to \Z_p$, the
{\em twisted Alexander matrix of $G$} is the $pm\times pn$ matrix
\[
A_{\lambda}=\big(\bar\lambda\phi\partial_i(r_k) \big)^{\b}
\]
obtained from $\big(\bar\lambda\phi\partial_i(r_k) \big)$
by replacing each entry $e$ with $\b(e)$.
\end{definition}

\begin{prop}
\label{prop:invfact}
Let $G$ be a commutator-relators group, and
let $K=\ker(\lambda:G\to \Z_p)$.  The matrix $A_{\lambda}$
is a relation matrix for the group $H_1(K)\oplus\Z^{p-1}$.
\end{prop}
A proof can be found in \cite{Jo}.
The matrix $A_{\lambda}$ is equivalent (via row-and-column operations)
to a diagonal matrix, from which the rank and elementary divisors of
$H_1(K)$ can be read off.

\subsection{Nilpotent quotients}
\label{subsec:nilpquot}
We now apply the above procedure to a particular class
of groups: the nilpotent quotients $G/G_q$, $q\ge 3$, of a
commutator-relators group $G=\F(n)/R$.

Let $\lambda:G/G_q\to\Z_p$ be a non-trivial representation.
To describe explicitly the presentation
matrix $A_{\lambda}$ of  Proposition~\ref{prop:invfact},
we need to examine more closely the Fox derivatives of
the relators $c_h$ and $r_k^{(q)}$ in the
presentation \eqref{eq:gmodgq} for $G/G_{q}$.

If $c$ is a non-simple basic commutator, then
Lemma~\ref{lem:foxcalc}~(\ref{foxa}), (\ref{foxc})
gives $\bar\lambda\phi(\d_i c)=0$.   If $c=[x_{i_1}, x_{i_2},
\dots,x_{i_q}]$ is a simple commutator, then it follows from
Lemma~\ref{lem:foxcalc}~(\ref{foxc}) that $\bar\lambda\phi\left(\partial_i
c\right)$ is either zero or of the form $e=\pm (\zeta ^{a_1}-1)\cdots
(\zeta ^{a_{q-2}}-1)\in\Z\Z_p$, for some integers $1\le a_j\le p-1$.

Recall that the truncation $r_k^{(q)}$ is a product of
basic commutators of length $<q$. The same argument shows that
$\bar\lambda\phi\big(\d_i r_k^{(q)}\big)$ is a linear combination
of elements in $\Z\Z_p$ of the form
$(\zeta ^{a_{i_1}}-1)\cdots (\zeta ^{a_{i_j}}-1)$, for $j<q-2$.

The following lemma shows the typical simplifications that
we will perform on $\big(\bar\lambda\phi(\d c_h)\big)^{\b}$ and
$\big(\bar\lambda\phi(\d r_k^{(q)})\big)^{\b}$.

\begin{lem}
\label{lem:diag}
The integral $p\times p$ matrix $e^{\b}$ of
$e=(\zeta ^{a_1}-1)\cdots(\zeta ^{a_{k}}-1)\in\Z\Z_p$
has diagonal form
\[
\big(\underbrace{p^{r-1}, \dots, p^{r-1}}_{p-l-1},
\underbrace{p^r, \dots, p^r}_{l}, 0 \big),
\]
where $r=\lceil\frac{k-1}{p-1}\rceil$, and $l=k-1-(r-1)(p-1)$.
Moreover, there is a sequence of row and column operations,
independent of the particular $e$, that brings $e^{\b}$ to that
diagonal form. \qed\renewcommand{\qed}{}
\end{lem}

\begin{prop}
\label{prop:freenilp}
Let $K$ be an index~$p$ normal subgroup of $\F(n)/\F(n)_q$. Then:
\[
H_1(K) = \Z^n \oplus (\Z/p^{r-1}\Z)^{(n-1)(p-l-1)} \oplus
(\Z/p^{r}\Z)^{(n-1)l},
\]
where $r=\left\lceil{\tfrac{q-2}{p-1}}\right\rceil$, and $l=q-2-(r-1)(p-1)$.
\end{prop}

\begin{proof}
In this case, only commutator relators are present, so
Lemma~\ref{lem:diag}, applied to each entry $\bar\lambda\phi(\d c_h)$,
shows that the matrix $A_{\lambda}$ is equivalent
the following diagonal matrix:
\begin{equation}
\label{eq:Dmat}
\qquad\qquad\qquad\quad
D=\big(\underbrace{p^{r-1}, \dots, p^{r-1}}_{(n-1)(p-l-1)},
\underbrace{p^r, \dots, p^r}_{(n-1)l},
\underbrace{0,\dots,0}_{n+p-1} \big).
\qquad\qquad\quad\qed
\end{equation}
\renewcommand{\qed}{}
\end{proof}

\begin{thm}
\label{thm:h1sgp}
Let $G=\F(n)/R$ be a commutator-relators group.
Let $K$ be an index~$p$ normal subgroup of $G/G_q$.
Set $r=\left\lceil{\tfrac{q-2}{p-1}}\right\rceil$.  Then:
\[
H_1(K) = \Z^n \oplus \bigoplus_{i=0}^{r} (\Z/p^{i}\Z)^{d_{i}},
\]
for some positive integers $d_0, \dots, d_{r}$ such that
$d_0+ \cdots +d_{r}=(n-1)(p-1)$ and $d_r\le l(n-1)$.
\end{thm}

\begin{proof}
Let $K=\ker (\lambda:G/G_q\to\Z_p)$. Consider the relation matrix
$A_{\lambda}$, corresponding to the presentation $G/G_q=\F/R\F_q$
from \eqref{eq:gmodgq}.  Partition $A_{\lambda}$ into two blocks,
$A_{\lambda}=\bigl(\begin{smallmatrix}B_{\lambda} \\
C_{\lambda}\end{smallmatrix}\bigr)$, where $B_{\lambda}$
corresponds to the relators $R$, and $C_{\lambda}$ corresponds
to the basic commutators.

Assume that the row and column operations of Lemma~\ref{lem:diag} have
already been performed. Then, after moving all the zero columns to the
right, $A_{\lambda}$ is equivalent to
$\Big(\begin{smallmatrix}
B_{\lambda}' & 0 \\
D' & 0 \\
\end{smallmatrix}\Big)$, where $D=\left(D'\ 0\right)$ is the diagonal
matrix~(\ref{eq:Dmat}). Since the number of zero diagonal elements
of $D$ is $n+p-1$, the rank of $K$ is $n$.  Since the non-zero
diagonal elements of $D$ are either $p^{r-1}$ or $p^r$, the elementary
divisors of $K$ are among $p,p^2,\dots ,p^r$.
\end{proof}

\subsection{$\nu$-Invariants}
\label{subsec:nuinv}
In view of Theorem~\ref{thm:h1sgp}, we define the following numerical
invariants of isomorphism type for the nilpotent quotients of a group.

\begin{definition}
\label{defin:nuinvar}
Let $G$ be a commutator-relators group, and $G/G_q$ be the $(q-1)^{\text
st}$ nilpotent
quotient of $G$. Given a prime $p$, and a positive integer $d$, define
\[
\nu_{p,d}(G/G_q) = \#\big\lbrace K \lhd G/G_q \:\big|\: [G/G_q:K]=p
\text{ and } \dim_{\,\Z_p} (\Tors H_1(K))\otimes \Z_p = d \:\big\rbrace.
\]
\end{definition}

\begin{example}
If $q=3$, then $H_1(K) = \Z^n \oplus \Z_p^d$, for some $0\le d\le n-1$.
So we have invariants $\nu_{p,0}(G/G_3),\dots ,\nu_{p,n-1}(G/G_3)$
for the second nilpotent quotient of $G$.
Since $\sum_{d=0}^{n-1} \nu_{p,d}=\frac{p^{n}-1}{p-1}$,
it is enough to compute $\nu_{p,1},\dots ,\nu_{p,n-1}$.
\end{example}

\begin{example}
If $q=4$, and $p\ge 3$, then $H_1(K) = \Z^n \oplus \Z_p^d$, for some
$0\le d\le 2n-2$.   If $p=2$, then
$H_1(K)=\Z^n\oplus\Z_2^{d_1}\oplus\Z_4^{d_2}$,
for some $0\le d=d_1+d_2\le n-1$.
\end{example}

\subsection{Second nilpotent quotient}
\label{subsec:gmodg3}
We now restrict our attention to $G/G_3$.
From \eqref{eq:gmodgq}, for $q=3$ we obtain the presentation:
\begin{equation}
\label{eq:gmodg3}
G/G_3=\langle x_1, \dots,x_ n\mid r_1^{(3)},\dots, r_m^{(3)}, c_1,
\dots, c_l \rangle,
\end{equation}
where $l=2\binom{n+1}{3}$, and the basic commutators $c_h$ are of
the form $[x_i,[x_j,x_k]]$, with $j<k$ and $i\ge j$.

\begin{thm}
\label{thm:rankmodp}
Given an epimorphism $\lambda:G/G_3 \to \Z_p$, with kernel $K_{\lambda}$,
we have
\[
\dim_{\,\Z_p} (\Tors H_1(K_{\lambda}))\otimes \Z_p= n-1-\rank_{\Z_p}
M(\lambda).
\]
\end{thm}

\begin{proof}
Recall from the proof of Theorem~\ref{thm:h1sgp} that the relation
matrix of the abelian group $H_1(K_{\lambda})$ has the following form:
$A_{\lambda}=\Big(\begin{smallmatrix} B_{\lambda}' & 0 \\ C_{\lambda}' & 0 \\
\end{smallmatrix}\Big)$. We have already seen in
Proposition~\ref{prop:freenilp}
that $C_{\lambda}'$ is equivalent to a diagonal matrix
$D'=\Big(\begin{smallmatrix}
I_{(n-1)(p-2)} & 0
\\0 & p\cdot I_{n-1}
\end{smallmatrix}\Big)$.

Recall also that $r_k^{(3)}=\prod_{i<j}[x_i,x_j]^{\e_{i,j}(r_k)}$.
A computation using formula (\ref{foxa}) in  Lemma~\ref{lem:foxcalc} shows:
\begin{equation}\label{eq:rhotrunc}
\bar\lambda\phi\bigg(\frac{\partial r_k^{(3)}}{\partial x_l}\bigg)
=\sum_{i=1}^{n}\e_{i,l}(r_k)(\zeta ^{\lambda_i}-1),
\end{equation}
for $1\le l\le n$ and $1\le k\le m$.

Consider $e=\sum_{\s=1}^{p-1}a_{\s}(\zeta ^{\s}-1)\in \Z\Z_p$.
Set $a=\sum_{\s=1}^{p-1}a_{\s}$.
It is readily seen that the matrix $e^{\b}$ is equivalent to:
\begin{equation}
\label{eq:hbeta}
\begin{pmatrix}
 * & \cdots & * & p\cdot a  & 0 \\
\vdots & & \vdots&\vdots& \vdots \\
 * & \cdots & * &p\cdot a  & 0 \\
 * & \cdots & * &\sum_{\s=1}^{p-1}a_{\s}\cdot \s & 0
\end{pmatrix}.
\end{equation}

Now \eqref{eq:linalex}, together with \eqref{eq:rhotrunc} and
\eqref{eq:hbeta}, imply that $B_{\lambda}'$ is equivalent to
$\big(\begin{smallmatrix} * & 0 \\  * & M(\lambda)' \end{smallmatrix}\big)$,
where $M(\lambda)'$ is some codimension $1$ minor of $M(\lambda)$.
Hence, $A_{\lambda}$ is equivalent to:
\[
\begin{pmatrix}
* & 0 & 0 \\
* & M(\lambda)' & 0 \\
I_{(n-1)(p-2)} & 0 & 0 \\
0 & p\cdot I_{n-1} & 0
\end{pmatrix}.
\]

The theorem now follows from the following fact:
An integral matrix of the form
$\big(\begin{smallmatrix} Q \\ p\cdot I_n \end{smallmatrix}\big)$
is equivalent to
$\big(\begin{smallmatrix}I_r & 0 \\ 0 & p\cdot I_d \end{smallmatrix}\big)$,
where $r=\rank_{\Z_p} Q\otimes \id_{\Z_p}$, and $d=n-r$.
\end{proof}

\begin{cor}
\label{cor:nurho}
$\nu_{p,d}(G/G_3) = \#\{ K_{\lambda}\in \Sigma_p(G/G_3)\mid
\rank_{\Z_p} M(\lambda) =n-d-1\}$.
\end{cor}

\subsection{Resonance varieties and subgroups of $G/G_3$}
\label{subsec:rvsgp}
The following  theorem relates the distribution of
index~$p$ normal subgroups of $G/G_3$, according
to their abelianization, to the number of points on
the $n$-dimensional projective space over $\Z_p$,
according to the stratification by the resonance varieties.

\begin{thm}
\label{thm:nucount}
For $G$ a commutator-relators group with $H_2(G)$ torsion free,
\[
\nu_{p,d}(G/G_3)=\#(\PP_d(G,\Z_p)\setminus\PP_{d+1}(G,\Z_p)).
\]
\end{thm}

\begin{proof}
Follows from Corollary~\ref{cor:projres} and Corollary~\ref{cor:nurho}.
\end{proof}

\section{Complex arrangements}
\label{sec:cxa}
We illustrate the techniques developed in the previous sections
with the main example of spaces satisfying conditions~(\ref{cond1})
and (\ref{cond2}):  complements of complex hyperplane arrangements.

\subsection{Cohomology and fundamental group}
\label{subsec:linearr}
Let $\A'$ be a complex hyperplane arrangement, with complement $X'$.
Let $\A$ be a generic two-dimensional section of $\A'$, with complement $X$.  
Then, by the Lefschetz-type theorem of Hamm and L\^{e}~\cite{HL}, 
the inclusion $i:X\to X'$ induces an isomorphism 
$i_*:\pi_1(X)\to\pi_1(X')$ and a monomorphism   
$i^*:H^{2}(X')\to H^{2}(X)$.  By the Brieskorn-Orlik-Solomon 
theorem, the map $i^*$ is, in fact, an isomorphism. 
So, for our purposes here, we may restrict our attention to $\A$.

Let $\A =\{H_{1},\dots ,H_{n}\}$ be an arrangement of $n$ affine lines
in $\C^2$, in general position at infinity.  Let $v_1, \dots ,v_s$ be the
intersection points of the lines. If $v_q = H_{i_1}\cap \dots \cap H_{i_m}$,
set $V_q=\{i_1,\dots,i_{m}\}$ and $\overline{V}_q=V_q\setminus \{\max V_q\}$.
The level~$2$ of the lattice of $\A$  is encoded in the list
$L_2(\A)=\{V_{1},\dots,V_{s}\}$, which keeps track of the
incidence relations between the points and the lines of the
arrangement.

The following properties hold:

\begin{romenum}
\item \label{i}
The homology groups of $X=\C^2\setminus \bigcup_{i} H_i$ are
free abelian, of ranks $b_1=n$, $b_2=\sum_{q=1}^{s} \abs{\overline{V}_q}$,
and $b_{i}=0$ for $i>2$.  The cohomology ring is
determined by $L_2(\A)$ (see~\cite{OT}):
\[
H^{*}(X)  = \bigg( e_1,\dots ,e_n \: \bigg|\:
\begin{array}{l}%
e_i^2=0,\ e_ie_j=-e_je_i\\[0pt]
e_ie_j + e_je_k + e_ke_i = 0\ \text{ for } i,j,k\in V_q,\ 1\le q\le s
\end{array}
\bigg) .
\]

\item \label{ii}
The fundamental group $G=\pi_1(X)$ is a commutator-relators group:
\[
G=\langle x_1,\dots ,x_n \mid \b_q(x_i)x_i^{-1}=1 \quad
\text{for } i \in \overline{V}_q\, \text{ and }\, q=1,\dots ,s\rangle.
\]
The pure braid monodromy generators $\b_1,\dots,\b_s$ can be read off
from a `braided wiring diagram' associated to $\A$ (see \cite{CSbm}).
Moreover, the space $X$ is homotopy equivalent to the $2$-complex given
by this presentation (see \cite{Li1}).

\item \label{iii}
The second nilpotent quotient is determined by $L_2(\A)$:
\[
G/G_3 = \bigg\langle x_1,\dots ,x_n \: \bigg|\:
\begin{array}{ll}%
[x_i,\prod_{j\in V_q} x_j] &
\text{for } i \in \overline{V}_q,\quad 1\le q\le s  \\[0pt]
[x_i,[x_j,x_k]] &\text{for }  1\le j<k\le n, \quad j\le i\le n
\end{array}
\bigg\rangle .
\]
This follows from the presentation in (\ref{ii}), together with
\eqref{eq:gmodg3}
(see also~\cite{Ry1}).

\item \label{iv}
The linearized Alexander matrix is determined by $L_2(\A)$.
It is obtained by stacking $M_{V_1}(\lambda),
\dots, M_{V_s}(\lambda)$, where
$M_V(\lambda)$ is the $\abs{\overline{V}}\times n$ matrix with entries
\[
M_V(\lambda)_{i,j}= \delta_{j,V} \Big(\lambda_i -
\delta_{i,j} \sum_{k\in V} \lambda_k \Big),\qquad
\text{ for } i\in \overline{V} \text{ and } 1\le j \le n.
\]
For a detailed discussion of the Alexander matrix
and the Alexander invariant of $\A$, see \cite{CSai}.
\end{romenum}

From properties (\ref{i}) and (\ref{ii}), we deduce that $X$ satisfies the
conditions from Proposition~\ref{prop:twocw}.

\subsection{Resonance varieties over $\C$}
\label{subsec:cxrv}
The resonance varieties of a complex hyperplane arrangement
were introduced by Falk in~\cite{Fa}.
Let $\A$ be an arrangement of $n$ affine lines in $\C^2$,
in general position at infinity.  Set $\RR_{d}(\A):=\RR_d(X,\C)$.
By Theorem 3.1 in \cite{Fa}, this definition agrees with Falk's
definition.

Qualitative results as to the nature of the resonance varieties
of complex arrangements were obtained by a number of authors,
\cite{Yuz, Fa, CScv, Li2, LY}.
We summarize some of those results, as follows.

\begin{thm}
\label{thm:cxresvar}
Let $\RR_1(\A)\subset \C^n$ be the resonance variety of an
arrangement of $n$ complex hyperplanes.  Then:
\begin{alphenum}
\item \label{rva}
The ambient type of $\RR_1(\A)$ determines the isomorphism type of $H^{\le
2}(X)$.

\item \label{rvb}
$\RR_1(\A)$ is contained in the hyperplane $\Delta_n:=\{\sum_{i=1}^n\l_i=0\}$.
\item \label{rvc}
Each component $C_i$ of $\RR_1(\A)$ is a linear subspace.
\item \label{rvd}
$C_i\cap C_j=\{\bo\}$ for $i\ne j$.
\item \label{rve}
$\RR_{d}(\A)=\{\bo\}\cup\bigcup_{\dim C_i\ge d+1} C_i$.
\end{alphenum}
\end{thm}
\begin{proof}
Part~(\ref{rva}) was proved in \cite{Fa}.
Part~(\ref{rvb}) was proved in \cite{Yuz} and \cite{Fa}.
Part~(\ref{rvc}) was conjectured in \cite{Fa}, and proved in \cite{CScv}
and \cite{Li2}.
Part~(\ref{rvd}) is proved in \cite{LY}.
Part~(\ref{rve}) follows from \cite{LY}, Theorem~3.4,
as was pointed out to us by S.~Yuzvinsky.
\end{proof}

By Theorem~\ref{thm:reseqts}, the resonance varieties $\RR_d(\A)$
are the determinantal varieties associated to the linearized Alexander
matrix, $M$. For another set of explicit equations, obtained from a
presentation of the linearized Alexander invariant, see~\cite{CScv}.

All the components of $\RR_1(\A)$ arise from {\em neighborly partitions} of 
sub-arrangements of $\A$, see \cite{Fa}, \cite{LY}.  To a partition $\Pi$ of 
$\A'\subset \A$, such that a certain bilinear form associated to $\Pi$ is 
degenerate, there corresponds a component $C_{\Pi}$ of $\RR_1(\A)$.  
For each $V\in L_2(\A)$ with $\abs{V}\ge 3$, there is a {\em local} component,
$C_{V}=\Delta_n \cap \{\l_i=0\mid i\notin V\}$, of dimension $\abs{V}-1$, 
corresponding to the partition $(V)$ of $\A_V=\{H_i\mid i\in V\}$.  
The other components of $\RR_1(\A)$ are called {\em non-local}.  
For more details and examples, see \cite{Fa, CScv, Li2, LY}.

\subsection{Resonance varieties over $\Z_p$}
We now turn to the characteristic~$p$ resonance varieties,
$\RR_{d}(\A;\Z_p)$.  Recall that the variety $\RR_{d}(\A)$
has integral equations, so we may consider its reduction mod~$p$.
As we shall see, there are arrangements $\A$ such that
$\RR_{d}(\A;\Z_p)$ does not coincide with $\RR_{d}(\A)$ mod $p$,
for certain primes $p$.  Indeed:

\begin{itemize}
\item
The number of irreducible components, or the dimensions
of the components  may be different, as illustrated in
Examples~\ref{exm:maclane} and \ref{exm:monomial} below.

\item
The analogues of Theorem~\ref{thm:cxresvar} (\ref{rva}) and (\ref{rve})
fail in general, as seen in Examples~\ref{exm:diamond} and 
\ref{exm:monomial} below.
\end{itemize}
On the other hand, it seems likely that the analogues of
Theorem~\ref{thm:cxresvar} (\ref{rvb}), (\ref{rvc}) and (\ref{rvd}) hold for
every prime~$p$.

Now let $\nu_{p,d}(\A)=\nu_{p,d}(G/G_3)$ be the number of normal
subgroups of $G/G_3$ with abelianization $\Z^n\oplus \Z_p^d$,
for $0\le d\le n-1$.   By properties (\ref{i}) and (\ref{ii}) above,
Theorem~\ref{thm:nucount} applies, and so $\nu_{p,d}(\A)$ can
be computed from the $\Z_p$-resonance varieties.

\begin{cor}
\label{cor:rvmodp}
If $\RR_d(\A,\Z_p)=\RR_d(\A) \bmod p$, for all $d\ge 1$, then
\[
\nu_{p,d-1}(\A)=\frac{p^{d}-1}{p-1} m_{d},
\]
where $m_d$ is the number of components of $\RR_1(\A)$
of dimension $d$.
\end{cor}

\begin{proof}
From the assumption, properties (\ref{rvc})--(\ref{rve}) hold for
$\RR_d(\A,\Z_p)$. Therefore, $\PP_d(X,\Z_p)\setminus\PP_{d+1}(X,\Z_p)$
consists of $m_d$ disjoint, $d$-dimensional projective
subspaces in $\P(\Z_p^{n})$. The formula follows from
Theorem~\ref{thm:nucount}.
\end{proof}

In particular, if all the components of $\RR_1(\A)$
are local, the Corollary applies, with
$m_d=\#\{V\in L_2(\A)\mid |V|=d+1\}$.

\subsection{Examples}
\label{subsec:cxexm}
We conclude this section with a few examples that
illustrate the phenomena mentioned above. The motivation
to look at Examples~\ref{exm:diamond} and \ref{exm:monomial}
came from S. Yuzvinsky, who was the first to realize that there are
exceptional primes for these  arrangements.  His method of
computing the corresponding non-local components is different
from ours, though.

\begin{example}
\label{exm:braidarr}
Let $\A$ be the reflection arrangement of type $\operatorname{A}_3$,
with lattice
\[
L_2(\A)=\{123,145,246,356,16,25,34\}.
\]
The variety $\RR_1(\A)$ has $5$ components of dimension $2$.
The non-local component,
$C_{\Pi}=\{\l_1-\l_6=\l_2-\l_5=\l_3-\l_4=0\}\cap \Delta_6$,
corresponds to the partition
$\Pi=(1 6\mid 2 5\mid 3 4)$,
see \cite{Fa, CScv, Li2}.

For all primes $p$, Corollary~\ref{cor:rvmodp} applies, giving
$\nu_{p,1}=5(p+1)$.
\end{example}

\begin{example}
\label{exm:diamond}
Let $\A$ be the realization of the non-Fano plane, with lattice
\[
L_2(\A)=\{123,147,156,257,345,367,24,26,46\}.
\]
The variety $\RR_1(\A)$ has $9$ components of dimension $2$.
The non-local components are given by the partitions 
$\Pi_1=(1 3\mid 4 6\mid 5 7)$,
$\Pi_2=(1 5\mid 2 4\mid 3 7)$,
$\Pi_3=(1 7\mid 2 6\mid 3 5)$ 
of the corresponding type $\operatorname{A}_3$ sub-arrangements, 
see \cite{CScv}.

For $p>2$, Corollary~\ref{cor:rvmodp} applies, and so $\nu_{p,1}=9(p+1)$.

For $p=2$, though, $\RR_1(\A,\Z_2)$ has a single, $3$-dimensional
non-local component,
$C_{\Pi}=\{\l_1+\l_4+\l_7=\l_2+\l_5+\l_7=\l_3+\l_6+\l_7=0\}\cap\Delta_7$,
corresponding to $\Pi=(1\mid 3\mid 5\mid 7\mid 2 4 6)$.  Furthermore,
$\RR_2(\A,\Z_2)$ has a single, $1$-dimensional component, 
$C_{\Pi'}=\{\l_1+\l_7=\l_3+\l_7=\l_5+\l_7=\l_2=\l_4=\l_6=0\}$,
corresponding to $\Pi'=(1\mid 3\mid 5\mid 7)$, and $\RR_3(\A,\Z_2)=\{\bo\}$. 
Thus, $\nu_{2,1}=24$ and $\nu_{2,2}=1$.
\end{example}

\begin{example}
\label{exm:maclane}
Let $\A$ be one of the realizations of the MacLane matroid, with
\[
L_2(\A)=\{123,456,147,267,258,348,357,168,15,24,36,78\}.
\]
The variety $\RR_1(\A)$ has $8$ local components.  Despite the
fact that $\A$ supports many neighborly partitions,
$\RR_1(\A)$ has no non-local components, since Falk's degeneracy
condition is not satisfied, see~\cite{Fa}.

For $p\ne 3$, Corollary~\ref{cor:rvmodp} applies, and so $\nu_{p,1}=8(p+1)$.

For $p=3$, though, the degeneracy condition is satisfied, and
$\RR_1(\A,\Z_3)$ has a non-local, $2$-dimensional component,
$C_{\Pi}=\{\l_2+\l_5+\l_8=\l_3+\l_5-\l_8=\l_4-\l_5-\l_8=
\l_5-\l_6-\l_8=\l_1+\l_5=\l_7+\l_8=0\}$, corresponding to
$\Pi=(1 5 \mid 2 4 \mid 3 6 \mid 7 8 )$. Moreover, 
$\RR_2(\A,\Z_3)=\{\bo\}$.  Hence, $\nu_{3,1}=36$.
\end{example}

\begin{example}
\label{exm:monomial}
Let $\A$ be the realization of the affine plane over $\Z_3$, with lattice
\[
L_2(\A)=\{123,456,789,147,258,369,159,357,168,249,267,348\}.
\]
The variety $\RR_1(\A)$ has $12$ local components, and $4$ non-local
components of dimension~$2$, see \cite{Fa, CScv, Li2, LY}.

For $p\ne 3$, Corollary~\ref{cor:rvmodp} applies, and so $\nu_{p,1}=16(p+1)$.

On the other hand, $\RR_1(\A,\Z_3)$ has a single, $3$-dimensional
non-local component, 
$C_{\Pi}=\{\l_1+\l_6+\l_8=\l_2+\l_4+\l_9=\l_3+\l_5+\l_7=
\l_3+\l_4+\l_8=\l_3+\l_6+\l_9=\l_7+\l_8+\l_9=0\}$,
corresponding to $\Pi=(1 2 3 \mid 4 5 6\mid 7 8 9)$,  
or any other of the partitions that give rise to the 
$4$ non-local components of $\RR_1(\A)$. Moreover, 
$\RR_2(\A,\Z_3)=C_{\Pi}$, and $\RR_3(\A,\Z_3)=\{\bo\}$.  
Thus, $\nu_{3,1}=48$ and $\nu_{3,2}=13$.
\end{example}

\begin{example}
\label{exm:falkex}
Let $\A_1$ and $\A_2$ be generic plane sections of the two arrangements
from \cite{Fa}, Example~4.10.  Each arrangement consists of $7$
affine lines in $\C^2$, and each resonance variety has only local
components.  Thus, the $\nu$-invariants of $\A_1$ and $\A_2$ coincide.
On the other hand, as shown by Falk, there is no linear automorphism
$\C^7\to\C^7$ restricting to an isomorphism $\RR_1(\A_1)\to\RR_1(\A_2)$.
The same `polymatroid' argument shows that there is no automorphism
$\P(\Z_p^7)\to\P(\Z_p^7)$ restricting to
$\PP_1(\A_1,\Z_p)\to\PP_1(\A_2,\Z_p)$.
Thus, the ambient type of the (projective) resonance varieties carries more
information than the count of their points.
\end{example}

\section{Real arrangements}
\label{sec:real}
We conclude with an application to the classification of
arrangements of transverse planes in $\R^4$.  Though similar
in some respects to central line arrangements in $\C^2$,
such arrangements lack a complex structure.  That difference
manifests itself in the nature of the resonance varieties.

\subsection{Arrangements of real planes}
\label{subsec:planes}
A {\em $2$-arrangement} in $\R^{4}$ is a finite collection
$\A=\{H_1,\dots,H_n\}$ of transverse planes through the origin
of $\R^4$.  Such an arrangement $\A$ is a realization of the uniform
matroid $U_{2,n}$; thus, its intersection lattice is solely determined
by $n$.   Let $X=\R^{4}\setminus \bigcup_{i} H_i$ be the
complement of the arrangement.  The {\em link} of the
arrangement is $L=\SS^3\cap \bigcup_{i} H_i$.  Clearly,
the complement of $\A$ deform-retracts onto the complement
of $L$.  The link $L$ is the closure of a pure braid
in $P_n$, see \cite{Mz}, \cite{MS}.  Hence, $G=\pi_1(X)$
is a commutator-relators group.

The {\em linking numbers} of $\A$ are by definition
those of the link $L$.  They can be computed from
the defining equations of $\A$: If $H_i=\{\a_i=\a'_i=0\}$,
for some linear forms $\a_i, \a'_i:\R^4\to \R$,  then
$l_{i,j}=\sgn(\det(\a_i,\a'_i,\a_j,\a'_j))$, see~\cite{Zi}.
A presentation for the cohomology ring of $X$ in terms of the
linking numbers is given in~\eqref{eq:linkcohoring}.

Arrangements of transverse planes in $\R^4$ fall, in the terminology
of~\cite{DV}, into several types: horizontal and non-horizontal,
decomposable and indecomposable.
A $2$-arrangement $\A$ is {\em horizontal} if it admits a defining polynomial
of the form $f(z,w)=\prod_{i=1}^{n}(z+a_iw+b_i\bar{w})$.  From
the coefficients of $f$, one reads off a permutation
$\tau\in S_n$.  Conversely, given $\tau$, there is a horizontal
arrangement, $\A(\tau)$, whose associated permutation is $\tau$.
A $2$-arrangement is {\em decomposable} if its link
is the $(1,\pm 1)$-cable of the link of another $2$-arrangement,
and it is {\em completely decomposable} if its link can be obtained
from the unknot by successive $(1,\pm 1)$-cablings. See \cite{MS} for details.

\subsection{Resonance varieties}
\label{subsec:realrv}
Let $\RR_{d}(\A):=\RR_{d}(X,\C)$ be the $d^{\text{th}}$ resonance
variety of $\A$.  Recall that the resonance varieties form a tower
$\C^n=\RR_0\supset \RR_1\supset\cdots\supset \RR_{n-1}=\{\bo\}$.
Moreover, they are the determinantal varieties  of the $n\times (n-1)$
matrix $M(\l)$, whose entries are given by
$M(\l)_{k,j}=l_{k,j}\l_k-\delta_{k,j}(\sum_{i}l_{k,i}\l_i)$.

If $\A$ is decomposable, the top resonance variety,
$\RR_1(\A)$, contains as a component the hyperplane
$\D_n=\{\l_1+ \cdots +\l_n=0\}$.  Moreover, if $\A$ is completely
decomposable, $\RR_1(\A)$ is the union of a central arrangement
of $n-2$ hyperplanes in $\C^n$ (counting multiplicities), with
defining equations  of the form $\e_1\l_1+\cdots +\e_n\l_n=0$,
where $\e_i=\pm 1$.  If $\A$ is indecomposable, though, $\RR_1(\A)$
may contain non-linear components  (see Example~\ref{exm:indecomp}).

At the other extreme, all the components of the variety $\RR_{n-2}(\A)$ are
linear. It can be shown that a horizontal arrangement $\A$ is indecomposable
if and only if $\RR_{n-2}(\A)=\{\bo\}$.

\begin{example}
\label{exm:zieglerpair}
In \cite{Zi}, Ziegler provided the first examples of $2$-arrangements
with isomorphic intersection lattices, but non-isomorphic cohomology rings.
Those arrangements are: $\A=\A(1234)$ and $\A'=\A(2134)$.  We can  distinguish
their cohomology rings by counting the components of their resonance
varieties:
\[
\RR_1(\A)=\D_4 \, , \quad \RR_1(\A')=\D_4\cup \{\l_1+\l_2-\l_3-\l_4=0\}.
\]
The example $\A'$ shows that the analogues of
Theorem~\ref{thm:cxresvar}~(\ref{rvb}), (\ref{rvd}), (\ref{rve}) do not
hold for
$2$-arrangements:
\begin{itemize}
\item
The second component of $\RR_1(\A')$ does {\em not}
lie in the hyperplane  $\Delta_4$.

\item
The two components of $\RR_1(\A')$ do {\em not} intersect only at the origin,
but rather, in the $2$-dimensional subspace $\{\l_1+\l_2=\l_3+\l_4=0\}$.

\item
We have $\RR_2(\A')$
$=\{\l_1+\l_2=\l_3=\l_4=0\}\cup \{\l_1=\l_2=\l_3+\l_4=0\}$,
and thus the stratification of $\RR_1$ by $\RR_d$'s
is {\em not} by dimension  of components.
\end{itemize}
\end{example}

\begin{example}
\label{exm:decomp}
Let $\A=\A(321456)$ and $\A'=\A(213456)$.
Then:
\begin{align*}
\RR_1(\A)&=\RR_1(\A')=\Delta_6\cup \{\l_1+\l_2-\l_3-\l_4-\l_5-\l_6=0\},\\
\RR_2(\A)&=\RR_1(\A), \qquad \RR_2(\A')=\Delta_6.
\end{align*}
This example shows that the analogue of Theorem~\ref{thm:cxresvar}~(\ref{rva})
does not hold for $2$-arrangements:  The variety $\RR_1$ fails to determine
$\RR_2$, and thereby fails to determine the cohomology ring of the complement.
\end{example}

\begin{example}
\label{exm:indecomp}
The horizontal arrangement $\A(31425)$ is indecomposable.
Its resonance varieties are:
\begin{align*}
\RR_1=\:&\{
\begin{tabular}[t]{@{}p{3.87in}}
$\l_1^3
-\l_2^3
-\l_3^3
+\l_4^3
-\l_5^3
+\l_1^2\l_2
-\l_1\l_2^2
+\l_1^2\l_3
-\l_1\l_3^2
-\l_1^2\l_4
-\l_1\l_4^2
+\l_1^2\l_5
-\l_1\l_5^2
+\l_2^2\l_3
+\l_2\l_3^2
-\l_2^2\l_4
+\l_2\l_4^2
+\l_2^2\l_5
+\l_2\l_5^2
-\l_3^2\l_4
+\l_3\l_4^2
+\l_3^2\l_5
+\l_3\l_5^2
+\l_4^2\l_5
-\l_4\l_5^2
+2\l_1\l_2\l_3
-2\l_1\l_2\l_4
+2\l_1\l_2\l_5
-2\l_1\l_3\l_4
+2\l_1\l_3\l_5
-2\l_1\l_4\l_5
+2\l_2\l_3\l_4
-2\l_2\l_3\l_5
+2\l_2\l_4\l_5
+2\l_3\l_4\l_5
=0\}$
\end{tabular}
\\
\RR_2=\:
&\{\l_1+\l_2= \l_3= \l_4= \l_5=0\}\cup
\{\l_1+\l_3= \l_2= \l_4= \l_5=0\}\cup   \\
&\{\l_2+\l_4= \l_1=\l_3= \l_5=0\}\cup
\{\l_3+\l_4= \l_1=\l_2=\l_5=0\}\cup     \\
&\{\l_1+\l_5=\l_2= \l_3= \l_4= 0\}\cup
\{\l_4+\l_5= \l_1=\l_2= \l_3=0\}\cup  \\
&\{\l_1-\l_4= \l_2= \l_3= \l_5=0\}\cup
\{\l_2-\l_3= \l_1= \l_4= \l_5=0\}\cup   \\
&\{\l_2-\l_5=\l_1= \l_3= \l_4=0\}\cup
\{\l_3-\l_5=\l_1= \l_2= \l_4=0\}\\
\RR_3=\:&\{\bo\}
\end{align*}
This example shows that the analogue of Theorem~\ref{thm:cxresvar}~(\ref{rvc})
does not hold for $2$-arrangements:  The variety $\RR_1$ is not linear.
\end{example}

\subsection{Ziegler invariant}
\label{subsec:zinv}
The cohomology rings of the arrangements in Example~\ref{exm:zieglerpair}
were distinguished by Ziegler by means of an invariant closely related to
one of the $Z$-invariants introduced in \ref{subsec:handginv}.

Recall the sequence $0\to G_2/G_3\to G/G_3\to H\to 0$.
This central extension is determined by the map
$\chi^{\top}:G_2/G_3\to \L^2 H$, given explicitly by \eqref{eq:chitop}.
The invariant $Z_{0,1}(\chi)=\coker\chi^{\top}$ equals $H^2(G)=\Z^{n-1}$.
More information is carried by the next invariant,
\[
Z_{0,2}(\chi)=\coker\Big(\mu_{H}\circ\L^2\chi^{\top}: \L^{2} G_2/G_3
\to \L^{4}H\Big).
\]
Set $Z(\A):=Z_{0,2}(\chi)$.
It can be shown that $Z(\A)=\Z^{\binom{n-1}{3}-r}\oplus \Z_2^{r}$,
where $r$ is some integer that can be read off from the linking
graph $\mathcal{G}$ of the link of $\A$.

For example, $Z(\A(1234))=\Z$ and $Z(\A(2134))=\Z_2$, showing again
that the two arrangements have different cohomology rings.
But $Z(\A)$ is not a complete invariant of the cohomology ring.
For example, $Z(\A(21435))=Z(\A(31425))=\Z_2^4$, although the
two arrangements are distinguished by the $\nu$-invariants (see below).

\subsection{Classification for $n\le 6$}
\label{subsec:classify}
Let $G$ the group of an arrangement of $n$ transverse
planes in $\R^4$, and $G/G_3$ its second nilpotent quotient.
As can be seen in Table~\ref{tab:nuinv},
the $\nu_{3,d}$-invariants completely classify the
second nilpotent quotients (and, thereby the cohomology rings)
of $2$-arrangement groups, for $n\le 6$, with a lone exception.

\begin{table}
\[
\begin{array}
{|c||c|c|c|c|c|c|}
\hline
\multicolumn{1}{|c||}{n}
&\multicolumn{1}{|c|}{\A}
&\multicolumn{1}{|c|}{\nu_{3,0}}
&\multicolumn{1}{|c|}{\nu_{3,1}}
&\multicolumn{1}{|c|}{\nu_{3,2}}
&\multicolumn{1}{|c|}{\nu_{3,3}}
&\multicolumn{1}{|c|}{\nu_{3,4}}
\\
\hline\hline
3 & \A(123)   &  9 &  4 & \multicolumn{2}{c}{} & \\
\hline
4 & \A(1234)  & 27 &  0 & 13 &\multicolumn{1}{c}{} & \\
  & \A(2134)  & 18 & 20 &  2 &\multicolumn{1}{c}{} & \\
\hline
5 & \A(12345) & 81 &  0 &  0 & 40 & \\
  & \A(21345) & 54 & 27 & 35 &  5 & \\
  & \A(21435) & 36 & 66 & 17 &  2 & \\
  & \A(31425) & 51 & 60 & 10 &  0 & \\
\hline
6 & \A(123456) & 243 &    0 &    0 &   0 & 121 \\
  & \A(213456) & 162 &   81 &    0 & 107 &  14 \\
  & \A(321456) & 162 &    0 &  162 &  32 &   8 \\
  & \A(215436) & 108 &  126 &   87 &  38 &   5 \\
  & \A(214356) & 108 &  108 &  121 &  24 &   3 \\
  & \A(312546) &  72 &  186 &   90 &  14 &   2 \\
  & \A(341256) &  81 &  162 &  112 &   6 &   3 \\
		& \A(314256) & 117 &  162 &   74 &  10 &   1 \\
		& \A(241536) & 108 &  200 &   48 &   8 &   0 \\
  & \LL        &  81 &  162 &  112 &   6 &   3 \\
		& \MM        & 144 &  160 &   60 &   0 &   0 \\
\hline
\end{array}
\]
\caption{Arrangements of $n\le 6$ planes in $\R^4$:
Number $\nu_{3,d}$ of index~$3$ subgroups, according to their abelianization,
$\Z^n\oplus \Z_3^d$.}
\label{tab:nuinv}
\end{table}

The exception is Mazurovski\u{\i}'s pair, $\KK=\A(341256)$ and $\LL$.
The corresponding configurations of skew lines in $\R^3$ were introduced
in \cite{Mz}.  Explicit equations for $\KK$ and $\LL$ can be found
in \cite{MS}.  As noted in~\cite{Mz}, the links of $\KK$ and $\LL$
have the same linking numbers.   Thus,
$H^*(X_{\KK};\Z)\cong H^*(X_{\LL};\Z)$, and
$G_{\KK}/(G_{\KK})_{3}\cong G_{\LL}/(G_{\LL})_{3}$.
On the other hand, $G_{\KK}/(G_{\KK})_{4}\not\cong G_{\LL}/(G_{\LL})_{4}$,
as can be seen from the distribution of the abelianization of their index~$3$
subgroups, shown in Table~\ref{tab:index3sgp}.
\begin{table}
\[
\begin{array}
{|c||c|c|c|c|c|c|c|c|c|}
\hline
\multicolumn{1}{|c||}{d}
&\multicolumn{1}{|c|}{0}
&\multicolumn{1}{|c|}{1}
&\multicolumn{1}{|c|}{2}
&\multicolumn{1}{|c|}{3}
&\multicolumn{1}{|c|}{4}
&\multicolumn{1}{|c|}{5}
&\multicolumn{1}{|c|}{6}
&\multicolumn{1}{|c|}{7}
&\multicolumn{1}{|c|}{8}
\\
\hline\hline
G_{\KK}/(G_{\KK})_{4} & 81 & 0 & 162 & 0 & 112 & 0 & 6 & 0 & 3 \\
G_{\LL}/(G_{\LL})_{4} & 81 & 0 & 172 & 24 & 78 & 6 & 0 & 0 & 3 \\
\hline
\end{array}
\]
\caption{The groups $G_{\KK}/(G_{\KK})_{4}$ and $G_{\LL}/(G_{\LL})_{4}$:
Number of index~$3$ subgroups, according to their abelianization,
$\Z^6 \oplus \Z_3^d$.}
\label{tab:index3sgp}
\end{table}

We summarize the above discussion, as follows:

\begin{thm}
\label{thm:class6}
Let $(\A,\A')\ne (\KK,\LL)$ be a pair of $2$-arrangements of
$n\le 6$ planes in $\R^4$.  Then  $H^*(X)\cong H^*(X')$
if and only if $X\simeq X'$.
\end{thm}

In other words, up to $6$ planes, and with the exception of
Mazurovski\u{\i}'s pair, the classification of complements
of $2$-arrangements up to cohomology-ring isomorphism coincides
with the homotopy-type classification.  As shown in \cite{MS}, the latter
coincides with the isotopy-type classification, modulo mirror images.


\end{document}